\documentclass{amsart}
\newtheorem{theorem}{Theorem}[section]
\newtheorem{lemma}[theorem]{Lemma}
\newtheorem{remark}[theorem]{Remark}
\newtheorem{corollary}[theorem]{Corollary}

\theoremstyle{definition}

\numberwithin{equation}{section}

\def\t{\tau}
\def\il{\int\limits}
\def\g{\gamma}
\def\f{\frac}
\def\q{\quad}

\def\r{\rho}
\def\P{\Phi}
\def\l{\lambda}

\def\s{\sigma}
\def\ep{\varepsilon}
\def\d{\delta}
\def\ba{\begin{array}}
\def\pr{\prime}

\begin{document}

\begin{center}
{\Large\bf Continuous modified Newton's-type method for}\\
{\Large\bf nonlinear operator equations}
\end{center}

\vspace{0.5cm}
\begin{center}
%    Information for the first and second authors
\hskip 5mm Alexander G. Ramm \footnote{This paper was finished when
AGR was visiting
Institute for Theoretical Physics, University of Giessen. The author
thanks DAAD for support}  \hskip 2.5cm Alexandra B. Smirnova\\
E-mail: ramm@math.ksu.edu \hskip 1.5cm E-mail: smirn@cs.gsu.edu\\
\hskip 6mm Department of Mathematics \hskip 1.5cm Department of Math and Stat\\
\hskip 3mm Kansas State University \hskip 2.2cm Georgia State
University\\\hskip 1mm Manhattan, KS 66506, USA
\hskip 1.8cm Atlanta, GA 30303, USA\\
\vspace{0.3cm}

%    Information for the third  author
Angelo Favini\\
E-mail: favini@dm.unibo.it\\
Dipartimento di Matematica\\
Universita di Bologna\\540127 Bologna, Italy\\
\end{center}
\vspace{5mm}

\begin{center}
{\large\bf Abstract}\vspace{5mm}
\end{center}

\noindent {\it A nonlinear operator equation $F(x)=0$, $F:H\to H,$
in a Hilbert space  is considered. Continuous Newton's-type
procedures based on a construction of a dynamical system with the
trajectory starting at some initial point $x_0$ and becoming
asymptotically close to a solution of $F(x)=0$  as $t\to +\infty$
are discussed. Well-posed and ill-posed problems are investigated.
} \vspace{5mm}

\noindent {\bf Key words:} nonlinear problem, integral inequality,
Fr\'echet derivative, Newton method.

\noindent {\bf AMS subject classification:} 65J15, 58C15,
47H17\vspace{0.3cm}

\section{{\bf Introduction}}
\setcounter{equation}{0} \setcounter{theorem}{0}
\renewcommand{\thetheorem}{1.\arabic{theorem}}
\renewcommand{\theequation}{1.\arabic{equation}}
\vspace*{-0.5pt}

The theme of this paper is solving nonlinear operator equations
of the form:
\begin{equation}\label{1.1}
F(x)=0,\q F: H\to H,
\end{equation}
in a real Hilbert space $H$. We consider a real Hilbert space for
the sake of simplicity: numerical algorithms for solving
(\ref{1.1}) in a complex Hilbert space can be treated similarly.
In order to approximate a solution to equation (\ref{1.1}) we use
the idea developed in \cite{ars1}-\cite{ar}, which consists of
constructing a dynamical system with the trajectory starting at
some initial point $x_0$ and converging to a solution of
(\ref{1.1}) as $t\to +\infty$. This idea, in its simplest form
goes back to A.Cauchy (steepest descent) and was proposed in \cite{c} for
solving
some optimization problems by a continuous analog of the gradient
method. In \cite{AR} a wide class of linear ill-posed problems
was studied by the dynamical systems method.
In \cite{g} a continuous Newton's scheme
\begin{equation}\label{1.2}
\dot x(t)=-[F'(x(t))]^{-1}F(x(t)),\q x(0)=x_0 \in H,
\end{equation}
was studied and a theorem establishing convergence with the
exponential rate was proved. A modified continuous Newton's method
is proposed in \cite{a}:
\begin{equation}\label{a1}
\dot x(t)=-J(t)F(x(t)),\q x(0)=x_0 \in H, \q J(0)\in L(H),\hskip
80 pt
\end{equation}
\begin{equation}\label{a2}
\dot J(t)=-\mu \bigl[F^{\pr*}(x(t))F'(x(t))J(t)+
J(t)F'(x(t))F^{\pr*}(x(t))\bigr]+2\mu F^{\pr*}(x(t)),
\end{equation}
where $\mu$ is a positive constant. System (\ref{a1})-(\ref{a2})
avoids the inversion of the Fr\'echet derivative $F'(x)$, which is
numerically difficult in some applications.

The regularized Gauss-Newton's-type algorithm with simultaneous
updates of the operator $[F^{\pr *}(x(t))F'(x(t))+\ep(t)I]^{-1}$
was proposed in \cite{rs}:
\begin{equation}\label{rs1}
\dot x(t)=-D(t)\bigl[F^{\pr*}(x(t))F(x(t))+\ep(t)(x(t)-x_0)\bigr],
\end{equation}
\begin{equation}\label{rs2}
\dot D(t)=-\bigl[(F^{\pr*}(x(t))F'(x(t))+\ep(t)
I)D(t)-I\bigr],\hskip 13pt
\end{equation}
$$
x(0)=x_0\in H,\q D(0)\in L(H),\q 0<\ep(t)\to 0\q \mbox{as}\q t\to
+\infty.
$$
It is shown that $x(t)$ converges to a solution of (\ref{1.1}) at
the rate $O(\ep(t))$. The convergence theorem is proved without
assuming monotonicity of $F$ and bounded
invertibility of $F'(x)$.

In \cite{ars2} and \cite{ar} a fairly general approach to the
analysis of continuous procedures in a Hilbert space was
developed. According to this approach one investigates a solution
to the Cauchy problem for a nonlinear operator-differential
equation by using differential inequalities. In the well-posed
case (the Fr\'echet derivative operator $F'(x)$ is boundedly
invertible in a ball, which contains one of the solutions) one
investigates the Cauchy problem for an autonomous equation:
\begin{equation}\label{1.3}
\dot x(t)=\P(x(t)),\q x(0)=x_0.
\end{equation}
The choice of $\P:H\to H$ yields a corresponding continuous
process. In the ill-posed case ($F'(x)$ has a nontrivial
null-space at the solution (\ref{1.1})
or is not boundedly invertible) a regularized continuous
procedure is required. For this reason the Cauchy problem for
the following equation is to be analyzed:
\begin{equation}\label{1.4*}
\dot x(t)=\P(x(t),t),\q x(0)=x_0,
\end{equation}
with $\P:H\times[0,+\infty)\to H$. If one takes
$$
\P(h,t):=-[F'(h)+\ep(t)]^{-1}(F(h)+\ep(t)(h-x_0)),
$$
then one arrives at a continuously regularized Newton's scheme
(CRNS). The convergence analysis of CRNS is done in \cite{ars2}
under the assumption that $F'(x) \geq 0$ as an operator in $H$. For
$$
\P(h,t):=-[F^{\pr*}(h)F'(h)+\ep(t)]^{-1}(F^{\pr*}(h)F(h)+\ep(t)(h-x_0))
$$
one obtains continuously regularized Gauss-Newton's scheme
(CRGNS). The convergence theorems for CRGNS (see \cite{ars1} and
\cite{ars2}) do not use any assumption about the location of the
spectrum of $F'(x)$. The absence of such assumption is made
possible by source-type conditions.

In section 2 of our paper we study a continuous analog of a
modified  Newton's method:
\begin{equation}\label{1.4}
\dot x(t)=-[F'(x_0)]^{-1}F(x(t)),\q x(0)=x_0 \in H,
\end{equation}
for solving well-posed nonlinear operator equation (\ref{1.1}).
Theorem~\ref{th1} establishes exponential convergence of
(\ref{1.4}) to a solution of  (\ref{1.1}). Process (\ref{1.4}) can be
used in practical
computations when calculating and inverting of $F'(x)$ at each
moment of time require a considerable effort. Another continuous
algorithm, investigated in section 2:
\begin{equation}\label{1.5}
\dot x(t)=-B(t)F(x(t)),\q x(0)=x_0 \in H, \q B(0)=B_0\in L(H),
\end{equation}
\begin{equation}\label{1.6}
\hskip -75 pt\dot B(t)=-F^{\pr*}(x(t))F'(x(t))B(t)+
F^{\pr*}(x(t)),
\end{equation}
can also be recommended in the above situation. It allows one to
update $[F'(x)]^{-1}$ continuously for $t\in[0,+\infty)$ without
actual inversion of the Fr\'echet derivative. In Theorem~\ref{th2}
the exponential convergence of (\ref{1.5})-(\ref{1.6})
to a solution of (1.1) is proved.

For many important inverse problems of the  form (\ref{1.1}) the
operator $F'(x)$ is not boundedly invertible. For such problems
the regularized version of algorithm (\ref{1.4}) is suggested in
section 3
\begin{equation}\label{1.7}
\dot x(t)=-[F'(x_0)+\ep(t)I]^{-1}(F(x(t))+\ep(t)(x(t)-x_0)),\q
x(0)=x_0 \in H,\q \ep(t)>0.
\end{equation}
The convergence analysis of continuous regularized method
(\ref{1.7}) is done in Theorem~\ref{th3} under the following
assumption:
\begin{equation}\label{1.8}
F'(x)=F'(x_0)G(x_0,x),
\end{equation}
where
\begin{equation}\label{1.9}
||G(x_0,x)-I||\le C(G)||x_0-x||,\q x_0,x\in U(\r,\hat x),
\end{equation}
and $\hat x$ is a solution to (\ref{1.1}). Assumption
(\ref{1.8})-(\ref{1.9}) is similar to condition (8) in \cite{k}.
It means that the operator $F'(x)$ remains in principle the same
for all $x$ in a neighborhood of a solution up to some
modification by a linear operator $G(x_0,x)$. The reader may
consult \cite{k} for several examples of nonlinear inverse
problems for which condition (\ref{1.8})-(\ref{1.9}) can be
verified. As a consequence of Theorem~\ref{th3} we obtain the
stability of process (\ref{1.7}) towards noise in the data and
choose an optimal regularization parameter (the stopping time)
such that the method converges to a solution of (\ref{1.1}) when
the noise level tends to zero.
   
Our main motivations for this investigation are: 

1) We think that the dsm (dynamical systems method) that we 
develop in this paper (and in the earlier publications,cited in the 
references) is not only of theoretical interest, but also provide a 
powerful numerical tool for solving very wide variety of 
problems, namely all the problems which can be described by equation (1.1) 
with the nonlinearity satifying the assumptions of our theorems
formulated in Sections 2 and 3.

2) We think that the idea of constructing a method for solving
equation (1.1) which does not require inverting $F'(u)$
(see, for example, equations (2.21) and (2.22) below)
 is of both theoretical and practical interest even for well-posed 
problems. 

3) The dsm gives a general approach to constructing convergent iterative 
methods for solving ill-posed nonlinear problems. We do not address this 
part of the dsm in our present paper, but it has been addressed in detail
in [4]. 

Finally we note that the dsm was tested numerically (see [2], [3], [10]
for example), but it certainly of interest to study much more
the numerical performance of dsm. In this paper, however, 
the authors deal with the theoretical questions.

\section{{\bf Continuous modified Newton's schemes for well-posed problems}}
\setcounter{equation}{0} \setcounter{theorem}{0}
\renewcommand{\thetheorem}{2.\arabic{theorem}}
\renewcommand{\theequation}{2.\arabic{equation}}
\vspace*{-0.5pt}

In this section we solve nonlinear operator equation (\ref{1.1})
under the assumption that the Fr\'echet derivative of the operator
$F$ is boundedly invertible in a ball which contains one of the
solutions. Let $x_0$ be an initial approximation for a solution
to (\ref{1.1}) and $x(t)$ be a trajectory of the autonomous
dynamical system
\begin{equation}\label{ch1}
\dot x(t)=\P(x(t)),\q 0\le t<+\infty,\q x(0)=x_0.
\end{equation}
Lemma~\ref{lem1} below (see \cite{ar})
gives simple sufficient conditions
on nonlinear operators $F$ in (\ref{1.1}) and $\P$ in
(\ref{ch1}) which guarantee that:

(a) initial value problem (\ref{ch1}) is uniquely solvable for
all $t\in [0,+\infty)$;

(b) the solution $x(t)$ tends to one of the solutions of
(\ref{1.1}) as $t\to +\infty$.

\begin{lemma}
Let $H$ be a real Hilbert space, $F,\P: H\to H$.

Suppose that there exist some positive numbers $c_1$ and $c_2$ such
that $F$ and $\P$ are Fr\'echet differentiable in
$U(r,x_0):=\{x\in H,\,\,||x-x_0||\le r\}$,
$r:=\f{c_2||F(x_0)||}{c_1}$ and $\,\forall h\in U(r,x_0)$ the
following conditions hold
\begin{equation}\label{ch2}
(F'(h)\P(h)),F(h))\le -c_1||F(h)||^2,
\end{equation}
and
\begin{equation}\label{ch3}
||\P(h)||\le c_2||F(h)||.
\end{equation}
Then:

\noindent 1. there exists a global solution $x=x(t)$ to problem
(\ref{ch1}) in the ball $U(r,x_0)$;

\noindent 2. there exists
\begin{equation}\label{ch4}
\lim _{t\to+\infty}x(t)=\hat x,
\end{equation}
where $\hat x$ is a solution to (\ref{1.1}) in $U(r,x_0),$ and
\begin{equation}\label{ch5}
||x(t)-\hat x||\le r e^{-c_1 t},
\end{equation}
\begin{equation}\label{ch6}
||F(x(t))||\le ||F(x_0)||e^{-c_1 t}.
\end{equation}
\label{lem1}\end{lemma}

{\bf Proof.} From the Fr\'echet differentiability of $\P$ the
local existence of a solution to (\ref{ch1}) follows, and from
(\ref{ch1}) one gets:
\begin{equation}\label{ch7}
\f{d}{dt} \{F(x(t))\}=F'(x(t))\dot x(t)=F'(x(t))\P(x(t)).
\end{equation}
Let $\l(t):=F(x(t))$. Then
\begin{equation}\label{ch8}
\dot \l(t)=F'(x(t))\P(x(t)),\q \l(0)=F(x_0).
\end{equation}
At least for sufficiently small $t$, for which $x(t)\in U(r,x_0)$,
one can use estimate (\ref{ch2}) and get:
\begin{equation}\label{ch8'}
\f{1}{2}\f{d}{dt}||\l(t)||^2=(\dot
\l(t),\l(t))=(F'(x(t))\P(x(t)),F(x(t)))\le -c_1||\l(t)||^2.
\end{equation}
Thus, at least for sufficiently small $t\ge
0,$ one gets:
\begin{equation}\label{ch9}
||\l(t)||\le ||F(x_0)||e^{-c_1t}.
\end{equation}
For $0\le t_1\le t_2$ by (\ref{ch3}) one has
$$
||x(t_2)-x(t_1)||\le\left\Vert\int^{t_2}_{t_1}\dot
x(s)ds\right\Vert\le\int^{t_2}_{t_1}||\P(x(s))||ds\le
c_2\int^{t_2}_{t_1}||\l(s)||ds
$$
\begin{equation}\label{ch10}
\le\f{c_2||F(x_0)||}{c_1}\left(e^{-c_1t_1}-e^{-c_1t_2}\right)<\f{c_2||F(x_0)||}{c_1}e^{-c_1t_1}.
\end{equation}
Setting $t_1=0$ and $t_2=t$, one concludes from (\ref{ch10}) that
$x(t)\in U(r,x_0)$ with $r=\f{c_2||F(x_0)||}{c_1}$ whenever it is
defined. Therefore the standard argument yields existence and
uniqueness of a solution to (\ref{ch1}) on $[0,+\infty)$. Now let
in (\ref{ch10}) $t_1=t$ and $t_2\to +\infty$ . Then one gets
(\ref{ch5}), and the limit $\hat x$ in (\ref{ch4}) does exist due
to (\ref{ch10}). From (\ref{ch9}) one concludes that $\hat x$ is a
solution to (\ref{1.1}). Inequality (\ref{ch6}) follows from
(\ref{ch9}). Lemma 2.1 is proved. \qed

\begin{remark} Assumptions of Theorem~\ref{th1} do
not imply uniqueness of a solution to equation (\ref{1.1}). If
(\ref{1.1}) is not uniquely solvable then $x(t)$ converges to one
of its solutions in $U(r,x_0)$.

In Lemma 2.1 we have assumed that $c_1$ and $c_2$
are some known constants in the ball $U(r,x_0)$, and this assumption 
allowed us to
define $r$ explicitly in terms of the ratio $\frac {c_2}{c_1}$ and 
$||F(x_0)||$.
One may assume that $\frac {c_2}{c_1}$ is not a constant but a function  
of $r$, $\frac {c_2}{c_1}:=c(r).$ In this case ,
in order that the argument of Lemma 2.1 be valid,
one has to satisfy the inequality
$c(r)||F(x_0)||\leq r$.
For example, if $\frac {c(r)}{r}\to 0$ as $r\to \infty$, 
then there always exists an
$r>0$ such that the conclusion of Lemma 2.1 holds 
and the assumptions of this lemma
are satisfied in the ball  $U(r,x_0)$.
 \end{remark}

Now consider the following continuous modified Newton's scheme:
\begin{equation}\label{2.2}
\dot x(t)=-[F'(x_0)]^{-1}F(x(t)),\q x(0)=x_0 \in H.
\end{equation}
Theorem~\ref{th1} below establishes a relation between the
asymptotic behavior of a solution $x(t)$ to (\ref{2.2}) and
solutions to equation (\ref{1.1}). It is a consequence of
Lemma~\ref{lem1}.

\begin{theorem}
Let $H$ be a real Hilbert space, $F: H\to H$. Assume that $F$  is
Fr\'echet differentiable, its Fr\'echet derivative $F'$ is
Lipschitz-continuous:
\begin{equation}\label{t1}
||F'(x_1)-F'(x_2)||\le M_2||x_1-x_2||\q \forall x_1,x_2 \in
U(\tilde r,x_0),
\end{equation}
where
\begin{equation}\label{t3}
U(\tilde r,x_0):=\{x\in H,\,\,||x-x_0||\le \tilde r\},\q \tilde
r:=\f{1}{2M_2m_1},\q m_1:=||[F'(x_0)]^{-1}||,
\end{equation}
and
\begin{equation}\label{t2}
4M_2m_1^2||F(x_0)||\le 1.
\end{equation}
Then

\noindent 1. There exists a unique solution $x=x(t)$, $t\in
[0,+\infty)$, to problem (\ref{2.2}).

\noindent 2. $x(t)\in U(\tilde r,x_0)\q \forall t\in [0,+\infty)$,
and
\begin{equation}\label{t4}
\lim _{t\to+\infty}x(t)=\hat x,
\end{equation}
where $\hat x$ is a solution to (\ref{1.1}).

\noindent 3. The following estimates hold:
\begin{equation}\label{t5}
||x(t)-\hat x||\le 2m_1||F(x_0)||e^{-\f{t}{2}},
\end{equation}
\begin{equation}\label{t6}
||F(x(t))||\le ||F(x_0)||e^{-\f{t}{2}}.
\end{equation}
\label{th1}\end{theorem}

{\bf Proof.} Take
\begin{equation}\label{*}
\P(x(t)):=-[F'(x_0)]^{-1}F(x(t)).
\end{equation}
Then, under assumptions (\ref{t1}) and (\ref{t3}) of Theorem~\ref{th1},
one gets $\forall h\in U(\tilde r,x_0)$
$$
(F'(h)\P(h),F(h))=-(F'(h)[F'(x_0)]^{-1}F(h),F(h))=-||F(h)||^2
$$
$$
+(\{I-F'(h)[F'(x_0)]^{-1}\}F(h),F(h))
=-||F(h)||^2+(\{F'(x_0)-F'(h)\}[F'(x_0)]^{-1}F(h),F(h))
$$
\begin{equation}\label{ch11}
\le -||F(h)||^2+M_2m_1\tilde r||F(h)||^2=-\f{1}{2}||F(h)||^2.
\end{equation}
Also one has $||\P(h)||\le m_1||F(h)||$. Thus conditions
(\ref{ch2}) and (\ref{ch3}) of Lemma~\ref{lem1} hold for any $h\in
U(\tilde r,x_0)$ with $\P$ defined in (\ref{*}), $c_1=\f{1}{2}$
and $c_2=m_1$. Hence
$\,\,r:=\f{c_2||F(x_0)||}{c_1}=2m_1||F(x_0)||$. From (\ref{t3})
and (\ref{t2}) one has $\,\,2m_1||F(x_0)||\le
\f{1}{2M_2m_1}:=\tilde r$. Therefore (\ref{ch2}) and (\ref{ch3})
are satisfied on $U(r,x_0)$, $r:=\f{c_2||F(x_0)||}{c_1}$,
$c_1=\f{1}{2}$, $c_2=m_1$. Applying Lemma~\ref{lem1}, one
completes the proof. \qed

\begin{remark}

\noindent (a) Choosing $\P(h)=-[F'(h)]^{-1}F(h)$ one gets
Continuous Newton's method. In this case $c_1=1,$
$c_2=\mu_1:=\sup_{x\in U(r,x_0)}||[F'(x)]^{-1}||,$
and Lemma~\ref{lem1} yields the convergence theorem for Continuous
Newton's method \cite{g}.

\noindent (b) Choosing $\P(h)=-F(h)$, one gets a simple iteration
method, for which condition (\ref{ch2}) means strict monotonicity
of $F$:  $F' \geq c_1>0,$ and $c_2=1$.

\noindent (c) $\P(h)=-[F'(h)]^*F(h)$ corresponds to the gradient
method. 

Here $c_2=M_1:= \sup_{x\in U(r, x_0)}||F'(x)||$,
and $c_1=\mu_1^{-1}$.

\noindent (d) $\P(h)=-[F^{\pr*}(h)F'(h)]^{-1}F^{\pr*}(h)F(h)$
yields Continuous Gauss-Newton's scheme. Here $c_1=1$,
$\,c_2=\mu_1^2 M_1$, where $\mu_1$ is the same as in (a) above,
and $M_1$ is  the same as in (c) above.
\end{remark}

In order to avoid inversion of the Fr\'echet derivative $F'(x(t))$
even at the initial moment $t=0,$ one can consider the following
 algorithm, which is a Cauchy problem for a system of two equations:
\begin{equation}\label{2.7}
\dot x(t)=-B(t)F(x(t)),\q x(0)=x_0 \in H, \q B(0)=B_0\in L(H),
\end{equation}
\begin{equation}\label{2.8}
\hskip -75 pt\dot B(t)=-F^{\pr*}(x(t))F'(x(t))B(t)+
F^{\pr*}(x(t)).
\end{equation}

Equation (2.22) is similar to equation (1.3) in \cite{AR}.
To prove  Theorem~\ref{th2} below we  use the following lemma,
which is an operator-theoretical version of the  Gronwall
inequality:

\begin{lemma}
Let
\begin{equation}\label{l.5}
\frac{d V}{d t} + A(t)V(t)=G(t),\q V(0)=V_0,
\end{equation}
where $A(t),$ $G(t),$ $V(t)\in L(H),$ $L(H)$ is the set of linear
bounded operators on $H$, and $H$ is a real Hilbert space. If
there exists a scalar function $\zeta(t)>0, \,\,\zeta \in
L^1_{loc}(0, \infty),$ such that
\begin{equation}\label{l.6}
(A(t)h,h)\ge \zeta(t)||h||^2 \q \forall h\in H,
\end{equation}
then
\begin{equation}\label{l.7}
||V(t)||\le
e^{-\il^t_0\zeta(p)dp}\left[\il^t_0||G(s)||e^{\il^s_0\zeta(p)dp}\,ds+||V(0)||\right].
\end{equation}
\label{lem2}\end {lemma}

{\bf Proof.} (see \cite{rs}) Take any $h\in H.$ Since $H$ is a
real Hilbert space one has:
$$
\hskip -1.7cm \frac{1}{2} \frac{d}{dt}
||V(t)h||^2=\left(\frac{dV}{dt}h,V(t)h\right)
$$
$$
\hskip 4cm =-(A(t)V(t)h, V(t)h)+(G(t)h,V(t)h)
$$
\begin{equation}\label{l.8}
\hskip 117pt \le -\zeta(t)||V(t)h||^2+||G(t)||\,\,
||h||\,\,||V(t)h||
\end{equation}
Denote $v(t):=||V(t)h||$. Inequality (\ref{l.8}) implies
\begin{equation}\label{l.9}
v\dot v\le-\zeta(t)v^2+||G(t)||\,\,||h||\,v.
\end{equation}
Divide this inequality by the nonnegative $v$ and get a linear
first-order differential inequality from which one gets
(\ref{l.7}). Lemma~\ref{lem2} is proved. \qed

\begin{theorem}
Let $H$ be a real Hilbert space, $F: H\to H$. Assume that:

\noindent 1.  $ U(R,x_0):=\{x\in H,\,\,||x-x_0||\le
R\}$, the operator $F$  is twice Fr\'echet differentiable, $F'(x)$
is boundedly invertible, and
\begin{equation}\label{tt1}
||F'(x)||\le M_1,\q||F''(x)||\le M_2,\q
||[F'(x)]^{-1}||^2\le\f{1}{c}\q \forall x \in U(R,x_0),
\end{equation}
where
\begin{equation}\label{tt3}
R:=\f{\g c}{2M_1M_2\s^2},\q  \g:=\f{1-||F'(x_0)B_0-I||}{2}>0,\q
\s:=\f{M_1}{c}+||B_0||.
\end{equation}

\noindent 2. Equation (\ref{1.1}) is solvable  in $U(R,x_0)$ (not
necessarily uniquely), and $\hat x$ is a solution.

\noindent 3. $F(x_0)$ satisfies the following condition
\begin{equation}\label{tt4}
\left\{\f{2M_1M_2}{c}\s^3||F(x_0)||\right\}^{1/2}\le \g.
\end{equation}
Then:

\noindent 1. there exists a unique solution $(x(t),B(t))$, $t\in
[0,+\infty)$, to problem (\ref{2.7})-(\ref{2.8});

\noindent 2. $x(t)\in U(R,x_0)\q \forall t\in
[0,+\infty)$;

\noindent 3. the following estimates hold
\begin{equation}\label{tt5}
||x(t)-\hat x||\le \f{\s||F(x_0)||}{\g}\,\,e^{-\g t},
\end{equation}
\begin{equation}\label{tt6}
||F(x(t))||\le ||F(x_0)||\,\,e^{-\g t},
\end{equation}
\begin{equation}\label{tt7}
||F'(x(t))B(t)-I||\le (M_2||F(x_0)|| \s^2t+||F'(x_0)B_0-I||) \,e^{-c
t},\q \mbox{if} \q c=\g,
\end{equation}
\begin{equation}\label{tt7a}
||F'(x(t))B(t)-I||\le
\left(\f{M_2
||F(x_0)||\s^2}{|c-\g|}+||F'(x_0)B_0-I||\right) \,e^{-\min\{\g,
c\} \,t},\, \q \mbox{if} \q c\neq\g.
\end{equation}
\label{th2}\end{theorem}

{\bf Proof.} Under the assumptions of theorem 2.6 there exists a
unique solution $(x(t),B(t))$ to (\ref{2.7})-(\ref{2.8}) on some
interval $[0,\tau]$, and, at least for sufficiently small $t>0,$
$x(t)\in U(R,x_0)$. Since $\forall x\in U(R,x_0)$ and $\forall
h\in H$
\begin{equation}\label{2.9}
(F'(x)F^{\pr *}(x)h,h)\ge \f{1}{||[F'(x)]^{-1}||^2}\,\,||h||^2\ge
c||h||^2,
\end{equation}
by (\ref{2.8}) and Lemma~\ref{lem2} one gets
$$
||B(t)||\le e^{-ct}\left [\int^t_0||F^{\pr
*}(x(s))||e^{cs}\,ds+||B(0)||\right].
$$
Thus by (\ref{tt1}) and (\ref{tt3})
\begin{equation}\label{2.10}
||B(t)||\le \f{M_1}{c}\left(1-e^{-ct}\right)+||B_0||e^{-ct}\le
\f{M_1}{c}+||B_0||:=\s.
\end{equation}
Let us analyze the initial value problem for  $w(t):=F(x(t))$. One
has
$$
\dot w(t)=F'(x(t))\dot x(t)=-F'(x(t))B(t)w(t).
$$
Therefore
\begin{equation}\label{2.11}
\dot w(t)+w(t)+[F'(x(t))B(t)-I]w(t)=0,\q w(0)=F(x_0).
\end{equation}
Denote $W(t):=F'(x(t))B(t)-I$. Then
$$
\dot W(t)=F''(x(t))\dot x(t)B(t)+F'(x(t))\dot B(t)
$$
$$
=-F''(x(t))B(t)F(x(t))B(t)+F'(x(t))[-F^{\pr*}(x(t))F'(x(t))B(t)+
F^{\pr*}(x(t))]
$$
$$=-F''(x(t))B(t)F(x(t))B(t)-F'(x(t))F^{\pr*}(x(t))W(t).$$
Consider the problem
\begin{equation}\label{2.12}
\dot W(t)+F'(x(t))F^{\pr*}(x(t))W(t)=-F''(x(t))B(t)F(x(t))B(t),
\end{equation}
\begin{equation}\label{2.13}
W(0)=F'(x_0)B(0)-I.
\end{equation}
From (\ref{2.9}), (\ref{2.12})-(\ref{2.13}) and Lemma~\ref{lem2}
one obtains the estimate
\begin{equation}\label{2.13'}
||W(t)||\le e^{-ct}\left
[\int^t_0||F''(x(s))B(s)F(x(s))B(s)||e^{cs}\,ds+||W(0)||\right].
\end{equation}
Assumptions 1 and 2 of Theorem~\ref{th2} yield
$$
||F(x(t))||\le ||F(x(t))-F(x_0)||+||F(x_0)-F(\hat x)||\le 2M_1R
$$
for all values of $t$ such that $x(t)\in U(R,x_0)$. 

Thus:
$$
||W(t)||\le \f{2M_1M_2
R}{c}\left(\f{M_1}{c}+||B(0)||\right)^2+||W(0)||=\f{2M_1M_2
R\s^2}{c}+||W(0)||=\gamma + ||W(0)||.
$$
Hence one gets from (\ref{tt3})
\begin{equation}\label{2.14}
||W(t)||\le \f{1-||W(0)||}{2}+||W(0)||=\f{1+||W(0)||}{2}.
\end{equation}
Now one has the following differential inequality
$$
\f{1}{2}\f{d}{dt}||w(t)||^2=-||w(t)||^2-(W(t)w(t),w(t))\le
-\f{1-||W(0)||}{2}||w(t)||^2=-\g ||w(t)||^2.
$$
Therefore
\begin{equation}\label{2.14'}
||w(t)||\le ||w(0)||\,\,e^{-\g t}
\end{equation}
for all values of $t$, such that $x(t)\in U(R,x_0)$. If $0\le
t_1\le t_2$, one obtains
$$
||x(t_2)-x(t_1)||\le \left\Vert\int^{t_2}_{t_1}\dot
x(s)\,ds\right\Vert\le
\left(\f{M_1}{c}+||B(0)||\right)\int^{t_2}_{t_1}||w(s)||\,ds
$$
\begin{equation}\label{2.15}
\le \f{\s ||F(x_0)||}{\g}\,\,\left(e^{-\g t_1} -e^{-\g
t_2}\right).
\end{equation}
From (\ref{2.15}) by conditions (\ref{tt3}) and (\ref{tt4}) one
gets:
\begin{equation}\label{2.16}
||x(t_2)-x(t_1)||\le \f{\g c }{2 M_1M_2\s^2}\,\,
\f{2M_1M_2\s^3||F(x_0)||}{c\g^2} \left(e^{-\g t_1} -e^{-\g
t_2}\right) \le R\,\left(e^{-\g t_1} -e^{-\g t_2}\right).
\end{equation}
Since $B(t)$ is bounded whenever it is defined, estimate
(\ref{2.16}) implies that there exists a unique solution
$(x(t),B(t))$ to (\ref{2.7})-(\ref{2.8}) on $[0,+\infty)$ and
$\forall t\in [0,+\infty)$ $x(t)\in U(R,x_0)$. Setting $t_1=t$ and
$t_2\to+\infty$ in (\ref{2.15}) one gets (\ref{tt5}).  Inequality
(\ref{tt6}) now follows from (\ref{2.14'}). Let us go back to
(\ref{2.13'}). By (\ref{2.14'}) one has
\begin{equation}\label{2.17}
||W(t)||\le e^{-ct}\left [\int^t_0||F''(x(s))|| \,
||B(s)||^2\,||F(x_0)||\,e^{(c-\g) s}\,ds+||W(0)||\right].
\end{equation}
Estimate (\ref{2.17}) implies (\ref{tt7})-(\ref{tt7a}). This
completes the proof. \qed

\section{{\bf Ill-posed case. Continuously regularized modified Newton's scheme}}
\setcounter{equation}{0} \setcounter{theorem}{0}
\renewcommand{\thetheorem}{3.\arabic{theorem}}
\renewcommand{\theequation}{3.\arabic{equation}}
\vspace*{-0.5pt}

In many important applications the Fr\'echet derivative operator
is not boundedly invertible, i.e. the problem is ill-posed. To
overcome this difficulty we suggest a regularized version of
algorithm (\ref{2.2}):
\begin{equation}\label{3.1}
\dot x(t)=-[F'(x_0)+\ep(t)I]^{-1}(F(x(t))+\ep(t)(x(t)-x_0)),\q
x(0)=x_0 \in H,\q 0<\ep(t),
\end{equation}
where $x_0$ is chosen so that $(F'(x_0)h,h)\ge 0$ $\forall h\in
H$. If such a choice is not possible for the original equation
$F(x)=0$, one may consider an auxiliary problem $\phi
(x):=F^{\pr*}(x_0)F(x)=0$. If $F$ is Fr\'echet differentiable, one
has $\phi'(x_0)=F^{\pr*}(x_0)F'(x_0)$ and $(\phi'(x_0))h,h)\ge 0$
$\forall h \in H.$ The last equation, in general, is not
equivalent to (\ref{1.1}). However every solution to (\ref{1.1})
solves $\phi(x)=0$.  The convergence analysis of
(\ref{3.1}) is done in the following theorem.

\begin{theorem}
Let $H$ be a real Hilbert space, $F: H\to H$,  equation
(\ref{1.1}) be solvable (not necessarily uniquely), and $\hat x$
be a solution to (1.1). Assume that:

\noindent 1. A positive function $\ep(t)\in C^1[0,+\infty)$
converges monotonically to zero as $t\to +\infty$, $\f{\dot
\ep(t)}{\ep(t)}$ is nondecreasing, and $\ep(0)>|\dot \ep(0)|$.

\noindent 2. $F$ is Fr\'echet differentiable, its
Fr\'echet derivative $F'$ is Lipschitz-continuous:
\begin{equation}\label{3.2}
||F'(x_1)-F'(x_2)||\le M_2||x_1-x_2||\q\mbox{and} \q
F'(x)=F'(x_0)G(x_0,x),
\end{equation}
where
\begin{equation}\label{3.3}
G(x_0,x)\in L(H),\q ||G(x_0,x)-I||\le C(G)||x_0-x||,\q \forall
x_1, \,\,x_2,\,\,x\in U(\r,\hat x),
\end{equation}
\begin{equation}\label{3.5}
U(\r,\hat x):=\{x\in H:\,||x-\hat x||\le \r\},\q C(G)>0,\q
\r:=\f{\ep(0)-|\dot \ep(0)|}{M_2+C(G)\ep(0)}.
\end{equation}

\noindent 3. $F'(x_0)$ is non-negative definite:
\begin{equation}\label{3.4}
(F'(x_0)h,h)\ge 0\q \forall h\in H,\q ||x_0-\hat x||<\r.
\end{equation}

\noindent 4. There exist $v\in H$ such that $\hat
x-x_0=F'(x_0)v$,
\begin{equation}\label{3.6}
\ep(0)-|\dot \ep(0)|\ge
[M_2+C(G)\ep(0)]\ep(0)\sqrt{\f{2||v||}{M_2}}.
\end{equation}

Then a unique solution $x=x(t)$ to problem (\ref{3.1}) exists for
all $t\in [0,+\infty)$ and
\begin{equation}\label{3.6'}
||x(t)-\hat x||\le \f{\ep(0)-|\dot
\ep(0)|}{\ep(0)[M_2+C(G)\ep(0)]}\ep(t).
\end{equation}
\label{th3}\end{theorem}

\begin{remark} Inequality (\ref{3.6}) can
always be satisfied if \\ $\sqrt{2||v||M_2}<1$. Indeed,
inequality (\ref{3.6}) is equivalent to
\begin{equation}\label{rr*}
1-\sqrt{2||v||M_2}\geq \frac {|\dot\ep(0)|}{\ep(0)} +
\ep(0) C(G)\sqrt {\f {2||v||}{M_2}}.
\end{equation}
Thus, if
%\begin{equation}\label{rr**}
%|\dot\ep(0)|\le\ep(0)\left\{1-\sqrt{2||v||M_2}-
%\ep(0) C(G)\sqrt{\f{2||v||}{M_2}}\right\},
%\end{equation}
%then inequality (\ref{3.6}) holds.
%If
$\sqrt{2||v||M_2} < 1,$ then inequality (\ref{rr*}) holds if
$\frac {|\dot\ep(0)|}{\ep(0)} $ and $\ep(0)$ are sufficiently
small. For $\ep(t)=a\,e^{-bt}$, $a,b>0$,
 inequality (\ref{3.6})  holds if the following inequality is
valid:
$$
b + a C(G)\sqrt{\f{2||v||}{M_2}}\le
1-\sqrt{2||v||M_2}.
$$
The foregoing inequality holds if $a$ and $b$ are positive and
sufficiently small and $1>\sqrt{2||v||M_2}$. Since a priori
$||v||$ is not known, in the numerical applications of the scheme
one has to try different functions $\ep(t)$ for (\ref{3.6}) to be
fulfilled.
\end{remark}

\begin{remark} Consider nonlinear integral equation of the first kind:
\begin{equation}\label{r1}
F(x):=\psi(x)-y=0,\q \psi(x)(t):=\int_0^1k(t,s)g(s,x(s))\,ds,\q
t\in [0,1],
\end{equation}
where $k(t,s)\in L^\infty((0,1)^2)$  and $g(s,u)$ is twice
continuously differentiable with respect to $u$ on $0\le s,t\le
1, \,\,-\infty<u<+\infty$. Suppose $F: \,H^1[0,1]\to L^2(0,1)$.
Then $$(F'(x)h)(t)=\int^1_0k(t,s)g_x(s,x(s))h(s)\,ds.$$ Introduce
the nonlinear operator $\phi(x):=F^{\pr*}(x_0)F(x)$, $\phi:
\,H^1[0,1]\to H^1[0,1]$, and solve the equation $\phi(x)=0$.
Clearly $\phi'(x_0)$ is non-negative definite, i.e. condition 3
of Theorem~\ref{th3} holds. Under the additional assumptions
$|g_x(s,x_0)|\ge \kappa
>0$ for any $s\in (0,1),$ and $g(s,u)\in
C^3((0,1)\times(-\infty,+\infty)),$
one can take
$\,\,(G(x_0,x)h)(s):=\f{g_x(s,x(s))}{g_x(s,x_0(s))}h(s)\,\,$ in
order to satisfy condition 2 of  Theorem~\ref{th3}. Indeed,
$$
(\phi'(x)h)(t)=(\phi'(x_0)G(x_0,x)h)(t),
$$
and for any $h\in H$ the following estimates are used in [7]:
$$
||(G(x_0,x)-I)h||_{L^2}=
\left\{\int_0^1\left[\f{\int_0^1g_{xx}(s,(x_0+\theta(x-x_0))(s))\,
d\theta(x(s)-x_0(s))h(s)}
{g_x(s,x_0(s))}\right]^2\,ds\right\}^{1/2}
$$
$$
\le\f{||g_{xx}||_{L^\infty}}{\kappa}||x-x_0||_{L^\infty}\,||h||_{L^2}.
$$
Also
$
\left\Vert\f{d}{ds}(G(x_0,x)-I)h\right\Vert_{L^2}
$
$$
=
\left\{\int_0^1\left[\left(\f{\int_0^1g_{sxx}(s,(x_0+\theta(x-x_0))(s))
+g_{xxx}(s,(x_0+\theta(x-x_0))(s))(x'_0+\theta(x'-x'_0))(s)d\theta}{g_x(s,x_0(s))}\right.\right.\right.
$$
$$
-\left.\f{\int_0^1g_{xx}(s,(x_0+\theta(x-x_0))(s))\, d\theta
(g_{sx}(s,x(s))+g_{xx}(s,x(s))x'(s))}
{g^2_x(s,x_0(s))}\right)(x(s)-x_0(s))h(s)
$$
$$
+\left.\left.\f{\int_0^1g_{xx}(s,(x_0+\theta(x-x_0))(s))\,
d\theta}
{g_x(s,x_0(s))}((x'(s)-x'_0(s))h(s)+(x(s)-x_0(s))h'(s))\right]^2\,ds\right\}^{1/2}
$$
$$
\le\f{||g_{sxx}||_{L^\infty}}{\kappa}||x-x_0||_{L^\infty}||h||_{L^2}
+\f{2||g_{xxx}||_{L^\infty}(||x'||_{L^2}+||x'_0||_{L^\infty})}{3\kappa}||x-x_0||_{L^\infty}||h||_{L^\infty}
$$
$$
+\f{||g_{xx}||_{L^\infty}||g_{sx}||_{L^\infty}}{\kappa^2}||x-x_0||_{L^\infty}||h||_{L^2}
+\f{||g_{xx}||^2_{L^\infty}||x'||_{L^2}}{\kappa^2}||x-x_0||_{L^\infty}||h||_{L^\infty}
$$
$$
\f{||g_{xx}||_{L^\infty}}{\kappa}(||x'-x'_0||_{L^2}||h||_{L^\infty}
+||x-x_0||_{L^\infty}||h'||_{L^2}).
$$
The $L^\infty(0,1)$-norms of $x-x_0$ and $h$ can be estimated by
their $H^1[0,1]$-norms times some constants, due to Sobolev's
embedding theorems.

Thus if one assumes that equation (\ref{r1}) is solvable, $\hat
x$ is its solution, and in a neighborhood of $\hat x$ there
exists $x_0$ such that
$$
\hat x-x_0=\phi'(x_0)v,\q \sqrt{2||v||M_2}<1,
$$
then a unique solution $x=x(t)$ to the problem
$$
\dot
x(t)=-[\phi'(x_0)+\ep(t)I]^{-1}(\phi(x(t))+\ep(t)(x(t)-x_0)),\q
x(0)=x_0 \in H,\q 0<\ep(t),
$$
exists for all $t\in [0,+\infty)$ and
$$
||x(t)-\hat x||=O(\ep(t)),
$$
provided that the above assumptions on $k(t,s)$ and $g(s,u)$ are
satisfied and the choice of $\ep(t)$ is made according to
(\ref{rr*}) with $\f{\dot\ep(t)}{\ep(t)}$ being nondecreasing.
\end{remark}

{\bf Proof of Theorem 3.1} First, from (\ref{3.4}) one
concludes that the operator $[F'(x_0)+\ep(t)I]^{-1}$ is bounded
$\forall t\ge 0$. Let us show that if $x=x(t)$ solves (\ref{3.1}),
then $x(t)\in U(\r,\hat x)$ with $\r$ introduced in (\ref{3.5}).
Assume the converse: there exists $T>0$ such that
\begin{equation}\label{3.7}
||x(t)-\hat x||<\r \q \forall t\in [0,T)\q\mbox{and}\q||x(T)-\hat
x||=\r.
\end{equation}
For any $t\in [0,T]$ one has
$$
\f{1}{2}\f{d}{dt}||x(t)-\hat
x||^2=-([F'(x_0)+\ep(t)I]^{-1}[F'(\hat x)(x(t)-\hat
x)+R_2(x(t),\hat x)+\ep(t)(x(t)-x_0)],x(t)-\hat x),
$$
where $||R_2(x(t),\hat x)||\le\f{M_2}{2}||x(t)-\hat x||^2.$ Thus
one gets
$$
\f{1}{2}\f{d}{dt}||x(t)-\hat x||^2\le-||x(t)-\hat
x||^2-([F'(x_0)+\ep(t)I]^{-1}(F'(\hat x)-F'(x_0))(x(t)-\hat
x),x(t)-\hat x)
$$
$$
-\ep(t)([F'(x_0)+\ep(t)I]^{-1}(\hat x-x_0),x(t)-\hat
x)+\f{M_2}{2\ep(t)}||x(t)-\hat x||^3.
$$
Condition 2 of Theorem~\ref{th3} and the estimate
$||[F'(x_0)+\ep(t)I]^{-1}F'(x_0)||\le 1$ yield
$$
\f{1}{2}\f{d}{dt}||x(t)-\hat x||^2\le-||x(t)-\hat
x||^2+||G(x_0,\hat x)-I||\,||x(t)-\hat x||^2
+\ep(t)||v||\,||x(t)-\hat x||
$$
$$
+\f{M_2}{2\ep(t)}||x(t)-\hat x||^3\le -(1-C(G)\r)||x(t)-\hat
x||^2+\ep(t)||v||\,||x(t)-\hat x||
$$
\begin{equation}\label{3.8}
+\f{M_2}{2\ep(t)}||x(t)-\hat x||^3.
\end{equation}
Introduce the notation
 $Q(t):=||x(t)-\hat x||$. Inequality
(\ref{3.8}) implies
\begin{equation}\label{3.9}
\dot Q(t)\le -(1-C(G)\r)Q(t)+\ep(t)||v||+\f{M_2}{2\ep(t)}Q^2(t),\q
Q(0)=||x_0-\hat x||.
\end{equation}
Take $f(t)=\f{Q(t)}{\ep(t)}$. By assumption 1 of Theorem~\ref{th3}
one obtains:
\begin{equation}\label{3.10}
\dot
f(t)\le-\left(1-C(G)\r-\f{|\dot\ep(0)|}{\ep(0)}\right)f(t)+||v||+\f{M_2}{2}f^2(t),\q
f(0)=\f{||x_0-\hat x||}{\ep(0)}.
\end{equation}
From (\ref{3.5}) and (\ref{3.10}) one concludes that
\begin{equation}\label{3.11}
\dot
f(t)\le-\f{M_2(\ep(0)-|\dot\ep(0)|)}{\ep(0)[M_2+C(G)\ep(0)]}f(t)+||v||+\f{M_2}{2}f^2(t),\q
f(0)=\f{||x_0-\hat x||}{\ep(0)}.
\end{equation}
Let
\begin{equation}\label{^}
C_1:=\f{M_2}{2},\q C_2:=\f{M_2(\ep(0)-|\dot
\ep(0)|)}{\ep(0)[M_2+C(G)\ep(0)]},\q C_3:=||v||.
\end{equation}
If $g(t)$ is a solution to the initial value problem
\begin{equation}\label{^*}
\dot g(t)=C_1g^2(t)-C_2g(t)+C_3,
\end{equation}
\begin{equation}\label{^**}
g(0)=f(0),
\end{equation}
then inequality (\ref{3.11}) yields
\begin{equation}\label{^***}
f(t)\le g(t),
\end{equation}
whenever $g(t)$ and $f(t)$ are both defined. By (\ref{3.5}),
(\ref{3.4}) and (\ref{^}) one has
\begin{equation}\label{^****}
f(0)=\f{||x_0-\hat x||}{\ep(0)}<\f{\r}{\ep(0)}=\f{\ep(0)-|\dot
\ep(0)|}{\ep(0)[M_2+C(G)\ep(0)]}=\f{C_2}{2C_1}.
\end{equation}
By (\ref{3.6}) the equation $C_1g^2-C_2g+C_3=0$ has at least one
real root. If there are two roots, the smaller root is a stable
equilibrium for problem (\ref{^*}), which implies $g(t)\le
g(0)=f(0)<\f{C_2}{2C_1}$. Otherwise $\tilde g:=\f{C_2}{2C_1}$ is
a solution to (\ref{^*}), and $g(t)<\f{C_2}{2C_1}$ since
$g(0)<\f{C_2}{2C_1}$. Therefore from (\ref{^***}) and
(\ref{^****}) one derives:
$$
f(t)\le g(t)<\f{\r}{\ep(0)}.
$$
Hence inequality (\ref{3.11}) and conditions (\ref{3.5}) and
(\ref{3.6}) yield:
\begin{equation}\label{3.12}
f(t)<
\f{\ep(0)-|\dot\ep(0)|}{\ep(0)[M_2+C(G)\ep(0)]}=\f{\r}{\ep(0)}.
\end{equation}
Thus
\begin{equation}\label{3.13}
||x(t)-\hat x||<\f{\r}{\ep(0)}\ep(t)\le \r\q \forall t\in [0,T],
\end{equation}
which contradicts (\ref{3.7}). Therefore $x(t)\in U(\hat x,\r)$
for any $t$, and by the standard argument one concludes that
$x(t)$ is defined on $[0,+\infty)$. Inequality (\ref{3.6'})
follows from (\ref{3.5}) and (\ref{3.13}). \qed

\begin{corollary}
In this corollary it is shown that if the data are noisy, then
the stopping time can be chosen so that the solution to the
Cauchy problem with noisy data approximates a solution to
(\ref{1.1}) stably, i.e. with the error going to zero as the
noise level goes to zero. Let the operator $F$ in (\ref{1.1})
have the following form
\begin{equation}\label{c1}
F(x):=\psi(x)-y.
\end{equation}
Assume that $\psi$ is given exactly and in place of $y$ we know a
$\d$-approximation $y_\d$, satisfying the inequality
\begin{equation}\label{c2}
||y-y_\d||\le \d.
\end{equation}
Then
$$
\f{1}{2}\f{d}{dt}||x(t)-\hat x||^2\le -(1-C(G)\r)||x(t)-\hat
x||^2+\left(\ep(t)||v||+\f{\d}{\ep(t)}\right)||x(t)-\hat x||
$$
\begin{equation}\label{c3}
+\f{M_2}{2\ep(t)}||x(t)-\hat x||^3.
\end{equation}
Take $\t_\d$ such that
$\ep(\t_\d)=\left(\f{\d}{||v||}\right)^{\f{1}{2}}$. For $t=\t_\d$
we get $\f{\d}{\ep^2(\t_\d)}=||v||$ and therefore $\forall t\in
[0,\t_\d]$
\begin{equation}\label{3.10d}
\dot
f(t)\le-\left(1-C(G)\r-\f{|\dot\ep(0)|}{\ep(0)}\right)f(t)+2||v||+\f{M_2}{2}f^2(t),\q
f(0)=\f{||x_0-\hat x||}{\ep(0)}.
\end{equation}
Thus one gets
\begin{equation}\label{c4}
||x(\t_\d)-\hat x||\le
\f{\r}{\ep(0)||v||^{\f{1}{2}}}\,\,\d^{\f{1}{2}},
\end{equation}
provided that conditions 1, 2, 3 of Theorem~\ref{th3} and
inequality
\begin{equation}\label{3.6d}
\ep(0)-|\dot \ep(0)|\ge
2[M_2+C(G)\ep(0)]\ep(0)\sqrt{\f{||v||}{M_2}}
\end{equation}
hold.
\end{corollary}

 \end{document}